\documentclass[11pt,reqno]{amsart}
\usepackage[T1]{fontenc}
\usepackage{lmodern,amsmath,amssymb,microtype}
\usepackage[hmargin=1.02in,vmargin=0.91in]{geometry}
\usepackage[hidelinks]{hyperref}
\newtheorem{theorem}{Theorem}[section]
\newtheorem{lemma}[theorem]{Lemma}
\newtheorem{proposition}[theorem]{Proposition}
\theoremstyle{remark}
\newtheorem*{remark}{Remark}
\numberwithin{equation}{section}
\newcommand{\R}{\mathbb{R}}
\newcommand{\dd}{\,\mathrm d}
\newcommand{\eps}{\varepsilon}
\allowdisplaybreaks[3]
\AtBeginDocument{%
  \addtolength{\abovedisplayshortskip}{0pt plus 1.5pt}%
  \addtolength{\belowdisplayshortskip}{0pt plus 1.5pt}%
}
\title[Time-decay rates toward a composite wave]{Time-Decay Rates toward a Rarefaction--Contact--Shock Composite Wave for the One-Dimensional Navier--Stokes--Fourier Equations}
\author{Le Chen}
\thanks{School of Mathematics, Statistics and Mechanics, Beijing University
of Technology, Beijing 100124, P. R. China
(\href{mailto:chenle@emails.bjut.edu.cn}{\texttt{chenle@emails.bjut.edu.cn}}).}
\date{}
\hypersetup{pdftitle={Time-Decay Rates toward a Rarefaction-Contact-Shock Composite Wave},
 pdfauthor={Le Chen},pdfsubject={One-dimensional compressible Navier-Stokes-Fourier equations},
 pdfkeywords={viscous shock wave, rarefaction wave, viscous contact wave, time decay}}
\begin{document}

\begin{abstract}
We consider the one-dimensional compressible Navier--Stokes--Fourier
equations when the end states are connected by a rarefaction wave,
a contact discontinuity, and a shock. For a fixed weak nonzero shock,
sufficiently weak rarefaction and contact waves, and small $H^2(\mathbb R)$
initial perturbations, we establish asymptotic stability with a
time-dependent shift on the viscous shock. The $L^4$ time average of
the uniform perturbation and the shift velocity over $[t,2t]$ is
bounded by $C_S(1+t)^{-\frac{1}{4}}$. We prove
pointwise-in-time bounds of order $(1+t)^{-\frac{1}{12}}$ on
the whole line, $(1+t)^{-\frac{1}{4}}$ near the shifted shock, and
$(1+t)^{-\frac{1}{6}}$ on $|x|\ge\eta t$ for each fixed $\eta>0$.
The proof uses energy estimates and comparison on finite time intervals.
\end{abstract}
\maketitle
\vspace{-0.6\baselineskip}
\begin{quote}

\end{quote}
\setlength{\jot}{2.5pt}
\enlargethispage{3pt}

\section{Introduction}

We consider the Cauchy problem for the one-dimensional compressible
Navier--Stokes--Fourier equations in Lagrangian coordinates,
\begin{equation}\label{eq:nsf}
 \begin{cases}
 v_t-u_x=0,\\
 u_t+p_x=\mu(\frac{u_x}{v})_x,\\
 E_t+(pu)_x=(\frac{\kappa\theta_x}{v}+\frac{\mu uu_x}{v})_x,
 \end{cases}
 \qquad x\in\R,\quad t>0.
\end{equation}
Here $v>0$, $u$, and $\theta>0$ denote the specific volume, velocity,
and absolute temperature, respectively. For a polytropic ideal gas,
the pressure and specific total energy are $p=\frac{R\theta}{v}$ and
$E=c_v\theta+\frac{u^2}{2}$, where $R>0$ is the gas constant and $c_v>0$ is the
specific heat at constant volume. The viscosity $\mu>0$ and heat
conductivity $\kappa>0$ are constants. We write $c_p=c_v+R$ and
$\gamma=\frac{c_p}{c_v}$.

We prescribe the initial data $q(0,x)=q_0(x)$, where $q=(v,u,\theta)$
and
\[
 \lim_{x\to\pm\infty}q_0(x)=q_\pm=(v_\pm,u_\pm,\theta_\pm),
 \qquad v_\pm,\theta_\pm>0.
\]
For strong solutions, subtracting $u$ times the momentum equation from
the energy equation gives
\begin{equation}\label{eq:temperature}
 c_v\theta_t+pu_x=\kappa(\frac{\theta_x}{v})_x+\frac{\mu u_x^2}{v}.
\end{equation}

We assume that the solution of the Euler Riemann problem with end
states $q_-$ and $q_+$ consists of a $1$-rarefaction wave from $q_-$
to $q_*$, a $2$-contact discontinuity from $q_*$ to $q^*$, and a
$3$-shock from $q^*$ to $q_+$. Across the contact discontinuity,
$u_*=u^*$ and $p_*=p^*$, where $p_*=p(v_*,\theta_*)$ and
$p^*=p(v^*,\theta^*)$. The shock states satisfy the Rankine--Hugoniot
relations, and its speed $\sigma$ satisfies the Lax condition
\[
 \frac{\sqrt{\gamma R\theta_+}}{v_+}<\sigma<
 \frac{\sqrt{\gamma R\theta^*}}{v^*}.
\]
We measure the wave strengths by
\[
 \delta_R=|v_*-v_-|,\qquad
 \delta_C=|v^*-v_*|,\qquad
 \delta_S=v_+-v^*>0.
\]
The rarefaction and contact curves give
\[
 \delta_R\le |q_*-q_-|\le C\delta_R,\qquad
 \delta_C\le |q^*-q_*|\le C\delta_C.
\]

Let $q^r(\frac{x}{t})$ be the self-similar Euler rarefaction wave and
$q^R(t,x)$ its smooth approximation. We use the viscous contact
profile $q^C(t,x)$ and the viscous shock profile $q^S(\xi)$ from
\cite[Section~2]{KVW}. With $X(t)$ defined by \eqref{eq:shift} and
$X(0)=0$, set
\begin{equation}\label{eq:profile}
 \bar q(t,x)=q^R(t,x)+q^C(t,x)+q^S(x-\sigma t-X(t))-q_*-q^*.
\end{equation}
For $t>0$, we write $\bar q^{\,r}(t,x)$ for the superposition in
\eqref{eq:profile} with $q^r(\frac{x}{t})$ in place of $q^R(t,x)$.

The stability of rarefaction waves was established for a
one-dimensional model of viscous gas in~\cite{MN} and for the
compressible Navier--Stokes equations in~\cite{LX}. For a weak
rarefaction wave of the isentropic system,
\cite[Theorem~2.4]{HHLX} gives decay of the perturbation in $L^4$
and of the relative Riemann invariants in $L^\infty$.
For viscous contact waves, asymptotic stability was proved in~\cite{HMX}.
Decay estimates were obtained in~\cite{YangContact,LWX}.
For viscous shocks in a heat-conducting gas, \cite{KM} treats small
zero-integral perturbations and \cite{LiuShock} describes the shock
translation and the outgoing diffusion waves. The pointwise Green
function method gives orbital decay for weak shocks with small
$L^1\cap H^3$ data~\cite{MZ}. With small $H^2$ data, \cite{single}
gives a uniform decay rate for a fixed weak viscous shock with a
time-dependent shift. A shift is also used in~\cite{HW} for the
vanishing-dissipation limit with single-shock data.

For composite waves, the stability of a viscous contact wave
combined with rarefaction waves was proved in~\cite{HLM}. The
full system with two viscous shocks was studied in~\cite{HM}.
Huang and Lee~\cite{HL} treat two viscous shocks with a viscous
contact wave between them. Kang and Vasseur~\cite{KV} proved
contraction up to a shift for barotropic shocks. The method was
used for a rarefaction and a shock in the barotropic
system~\cite{KVWbaro}. For the full system,
\cite[Theorem~1.1]{KVW} gives uniform convergence to a rarefaction,
a viscous contact wave and a shifted viscous shock.

Wang and Wang~\cite{WW} used a time-dependent shift to prove planar
shock stability for the three-dimensional barotropic system on
$\R\times\mathbb T^2$. Decay estimates for a planar rarefaction wave
and a planar shock in the full system on $\R^3$ were obtained
in~\cite{DHLX} and~\cite{HX}, respectively.

We first obtain a time-averaged bound. For the pointwise rates,
we compare solutions that agree on a spatial interval at time
$t-\ell$. An exponential weight bounds their difference up to
time $t$.
\begin{enumerate}
\item We use the first-order estimate in
      \cite[Proposition~3.2]{KVW}, keeping track of $\delta_R$ and
      $\delta_C$ for the fixed shock. Differentiating the residuals
      gives the second-order estimate. The time integral of the
      dissipation gives the averaged bound.
\item Near the shock, the rarefaction and contact tails are
      exponentially small. We can take $\ell$ proportional to $t$
      and use the single-shock bound \eqref{eq:single-input}.
      We subtract the two shift equations to estimate the
      difference in shock positions.
\item Away from the shock, we start the comparison from a
      constant state with a cutoff perturbation and use
      \eqref{eq:constant}. Over the cutoff interval the profile
      variation contributes at most
      $C\ell^{\frac{3}{2}}(1+t)^{-\frac{1}{2}}$ in $H^2$.
      On $|x|\ge\eta t$, the bound improves to
      $C_\eta\ell^{\frac{3}{2}}(1+t)^{-1}$. We can therefore
      take $\ell$ proportional to $(1+t)^{\frac{1}{3}}$ and
      $(1+t)^{\frac{2}{3}}$, respectively.
\end{enumerate}

\begin{theorem}\label{thm:decay}
Assume that the end states are connected by the rarefaction wave,
contact discontinuity, and shock described above. Fix a viscous
$3$-shock $q^S$ of sufficiently small nonzero strength $\delta_S$.
There exist positive constants $\eps_S$, $C_S$, and $\nu$ such that, if
\begin{equation}\label{eq:small}
 q_0-\bar q(0,\cdot)\in H^2(\R;\R^3),\qquad
 \|q_0-\bar q(0,\cdot)\|_{H^2}+(\delta_R+\delta_C)^{\frac{1}{4}}\le\eps_S,
\end{equation}
then \eqref{eq:nsf} admits a unique global strong solution in the small
perturbation class. Moreover, $q-\bar q\in C([0,\infty);H^2(\R))$,
and $v$ and $\theta$ are bounded away from zero uniformly in time
and space. For the shift $X(t)$ of \eqref{eq:shift},
\[
 \lim_{t\to\infty}\big(\|q(t)-\bar q(t)\|_\infty+|\dot X(t)|\big)=0.
\]
For $t\ge1$,
\[
 \left(\frac{1}{t}\int_t^{2t}
    \big(\|q(\tau)-\bar q(\tau)\|_\infty+|\dot X(\tau)|\big)^4\dd\tau
 \right)^{\frac{1}{4}}
 \le C_S(1+t)^{-\frac{1}{4}} .
\]
Moreover,
\begin{align}
 &\|q(t)-\bar q(t)\|_{L^\infty(\R)}
   \le C_S(1+t)^{-\frac{1}{12}},\label{eq:global-rate}\\
 &\sup\nolimits_{|x-\sigma t-X(t)|\le\nu t}|q(t,x)-\bar q(t,x)|+|\dot X(t)|
   \le C_S(1+t)^{-\frac{1}{4}},\label{eq:shock-rate}
\end{align}
for all $t\ge0$. For each $\eta>0$ there exists $C_{S,\eta}>0$
such that
\begin{equation}\label{eq:away-contact}
 \sup_{|x|\ge\eta t}|q(t,x)-\bar q(t,x)|
       \le C_{S,\eta}(1+t)^{-\frac{1}{6}},\qquad t\ge0.
\end{equation}
For $t\ge1$, \eqref{eq:global-rate}--\eqref{eq:away-contact} also hold
with $\bar q^{\,r}$ in place of $\bar q$, using the same shift $X(t)$.
The constants may depend on the fixed shock and the physical parameters.
The constant $C_{S,\eta}$ also depends on $\eta$.
\end{theorem}

\begin{remark}
The time average follows from the integrated dissipation. For
the pointwise estimates, the cutoff datum at time $t-\ell$ changes
with $t$. This is why the comparison uses the uniform bound
\eqref{eq:single-input} for all small $H^2$ perturbations of the
fixed shock.
\end{remark}

Section~\ref{sec:profiles} recalls the wave profiles and estimates
the residuals. Section~\ref{sec:energy} defines the shift and
derives the a priori estimates. Section~\ref{sec:comparison} gives
the comparison estimates, which are used in Section~\ref{sec:proof}
to prove Theorem~\ref{thm:decay}.

\section{Preliminaries}\label{sec:profiles}

We use standard norms over $\mathbb R$ throughout. Here, $C,c>0$
denote generic constants that may change from line to line, depending
on the fixed shock, physical parameters, and a compact set of positive
states.

\subsection{Approximate rarefaction wave}

Set $\lambda_1(q)=-\frac{\sqrt{\gamma R\theta}}{v}$, $w_-=\lambda_1(q_-)$,
and $w_*=\lambda_1(q_*)$. Let $w$ solve
\[
 w_t+ww_x=0,\qquad
 w(0,x)=w_0(x)=\frac{w_-+w_*}{2}+\frac{w_*-w_-}{2}\tanh x.
\]
Following~\cite[Section~2.1]{KVW}, we define $q^R(t,x)$ as the
state on the $1$-rarefaction curve through $q_-$ and $q_*$ with
$\lambda_1(q^R(t,x))=w(1+t,x)$. The profile $q^R$ satisfies the
Euler equations. If $\delta_R=0$, we take $q^R\equiv q_-=q_*$.
The first- and second-order estimates and the end-state bounds are
given in~\cite[Lemma~2.1]{KVW}.

The third derivative needed below follows from the characteristics.
For $x=y+(1+t)w_0(y)$, we have $w(1+t,x)=w_0(y)$ and
$\dd x=[1+(1+t)w_0'(y)]\dd y$. Thus
\begin{equation}\label{eq:burgers}
 \begin{gathered}
 w_x(1+t,x)=\frac{w_0'(y)}{1+(1+t)w_0'(y)},\qquad
 w_{xx}(1+t,x)=\frac{w_0''(y)}{[1+(1+t)w_0'(y)]^3},\\
 w_{xxx}(1+t,x)=\frac{w_0'''(y)}{[1+(1+t)w_0'(y)]^4}
     -\frac{3(1+t)[w_0''(y)]^2}{[1+(1+t)w_0'(y)]^5}.
 \end{gathered}
\end{equation}
Along the rarefaction curve $\theta v^{\gamma-1}$ is constant and
$\frac{\mathrm d\lambda_1}{\mathrm dv}
=-\frac{(\gamma+1)\lambda_1}{2v}>0$. Hence
$c\delta_R\le w_*-w_-\le C\delta_R$, and the state on the
rarefaction curve is a smooth function of its characteristic speed
with bounded derivatives through order three. With the derivatives of $w$
evaluated at $(1+t,x)$, \eqref{eq:burgers} gives
\begin{align*}
 |w_{xx}|+|w_{xxx}|
 &\le C\left(\frac{w_0'(y)}{[1+(1+t)w_0'(y)]^3}
        +\frac{(1+t)[w_0'(y)]^2}{[1+(1+t)w_0'(y)]^5}\right)
 \le \frac{Cw_0'(y)}{[1+(1+t)w_0'(y)]^3},
\end{align*}
where we used $|w_0''(y)|+|w_0'''(y)|\le Cw_0'(y)$. Since
$w_t=-ww_x$, the chain rule yields
\begin{equation}\label{eq:rare-infty}
 \begin{aligned}
 |q^R_t|+\sum_{j=1}^3|\partial_x^j q^R|
 &\le C(w_x+|w_{xx}|+|w_{xxx}|)
 \le \frac{Cw_0'(y)}{1+(1+t)w_0'(y)}
 \le C\min\{\delta_R,(1+t)^{-1}\}.
 \end{aligned}
\end{equation}

For $\delta_R>0$, we have $c\delta_R\le w_0'(y)\le C\delta_R$
on $|y|\le1$ and $|w_0''(y)|\ge2\tanh(1)w_0'(y)$ on $|y|>1$.
Changing variables from $y$ to $w_0'(y)$ on each tail and using
$\int_{\R} w_x\dd x=w_*-w_-$, we obtain
\begin{align}
 &\|q^R_{xx}\|_2^2+\|q^R_{xxx}\|_2^2+\|q^R_x\|_4^4\notag\\
 &\quad\le C\int_\R
        \frac{[w_0'(y)]^2}{[1+(1+t)w_0'(y)]^5}\dd y
       +C\delta_R\|w_x(1+t,\cdot)\|_\infty^3\notag\\
 &\quad\le C\left(
       \frac{\delta_R^2}{[1+(1+t)\delta_R]^5}
       +\int_0^{C\delta_R}\frac{z}{[1+(1+t)z]^3}\dd z
       +\frac{\delta_R^4}{[1+(1+t)\delta_R]^3}\right)
 \le \frac{C\delta_R^2}{(1+\delta_R t)^2}.
 \label{eq:rare-higher}
\end{align}
The primitive of the integrand on the last line is
$\frac{z^2}{2[1+(1+t)z]^2}$.

For $x\ge (1+t)w_*$, applying \eqref{eq:rare-infty} and the
end-state bound of~\cite{KVW} gives
\[
 |q^R-q_*|+\sum_{j=1}^3|\partial_x^j q^R|
 \le C\delta_R e^{-2[x-(1+t)w_*]}+Cw_0'(y)
 \le C\delta_R e^{-2[x-(1+t)w_*]},
\]
since $y\ge x-(1+t)w_*\ge0$. On $x\le(1+t)w_-$, the corresponding
bound holds with $q_-$ and $e^{-2[(1+t)w_--x]}$.

\subsection{Viscous contact wave and viscous shock wave}

We use the viscous contact profile of~\cite[Section~2.2]{KVW}.
Its temperature $\theta^C$ is self-similar in $\frac{x}{\sqrt{1+t}}$ and
connects $\theta_*$ to $\theta^*$. The profile satisfies
\[
 v^C=\frac{R\theta^C}{p_*},\qquad
 u^C=u_*+\frac{\kappa}{c_p}\frac{\theta^C_x}{\theta^C},\qquad
 \theta^C_t=\frac{\kappa p_*}{Rc_p}
             \left(\frac{\theta^C_x}{\theta^C}\right)_x .
\]
Thus $p^C=p_*$, $v^C_t=u^C_x$, and
\[
 c_v\theta_t^C+p^C u_x^C
 =(c_v+R)\theta_t^C=\kappa\left(\frac{\theta_x^C}{v^C}\right)_x .
\]
By~\cite[Lemma~2.2]{KVW},
\begin{equation}\label{eq:contact}
 |\partial_x^j(v^C,\theta^C)|+(1+t)^{\frac{1}{2}}|\partial_x^j u^C|
 \le C\delta_C(1+t)^{-\frac{j}{2}}e^{-\frac{cx^2}{1+t}},\qquad j\ge1 .
\end{equation}
On $x\le0$ and $x\ge0$, respectively, $|q^C-q_*|$ and
$|q^C-q^*|$ are bounded by $C\delta_C e^{-\frac{cx^2}{1+t}}$.
Differentiating the self-similar profile gives
\[
 (v^C,\theta^C)_t=-\frac{x}{2(1+t)}(v^C,\theta^C)_x,\qquad
 u_t^C=-\frac{u^C-u_*+xu_x^C}{2(1+t)},\qquad
 u_{tx}^C=-\frac{2u_x^C+xu_{xx}^C}{2(1+t)}.
\]
It follows from \eqref{eq:contact} and the formula for $u^C$ that
\begin{align*}
 &|(v^C,\theta^C)_t|+(1+t)^{\frac{1}{2}}|u_t^C|+(1+t)|u_{tx}^C|\\
 &\quad\le \frac{|x|\,|(v^C,\theta^C)_x|}{2(1+t)}
       +\frac{|u^C-u_*|+|x|\,|u_x^C|}{2(1+t)^{\frac{1}{2}}}
       +\frac{2|u_x^C|+|x|\,|u_{xx}^C|}{2}\\
 &\quad\le C\delta_C(1+t)^{-1}
      \left(1+\frac{|x|}{\sqrt{1+t}}\right)e^{-\frac{2cx^2}{1+t}}
 \le C\delta_C(1+t)^{-1}e^{-\frac{cx^2}{1+t}} .
\end{align*}
Consequently,
\[
 \|u_t^C(t)\|_{H^1}\le C\delta_C(1+t)^{-\frac{3}{2}}
   \left(\int_\R e^{-\frac{2cx^2}{1+t}}\dd x\right)^{\frac{1}{2}}
 \le C\delta_C(1+t)^{-\frac{5}{4}}.
\]

Let $q^S(\xi)=(v^S,u^S,\theta^S)(\xi)$, $\xi=x-\sigma t$,
be the viscous $3$-shock profile connecting $q^*$ to $q_+$.
With the translation fixed, we use the exponential bounds at
both ends and the estimates of~\cite[Lemma~2.3]{KVW}
\[
 0<(v^S)'\le C\delta_S^2e^{-c\delta_S|\xi|},\qquad
 |(q^S)'|+\delta_S^{-1}|(q^S)''|\le C(v^S)'.
\]
Integrating the traveling-wave equations from $q^*$ gives
$u^S-u^*=-\sigma(v^S-v^*)$. The integrated momentum equation
and the integrated energy equation with $u^S$ times the momentum
equation subtracted give
\begin{align*}
 -\mu\sigma\frac{(v^S)'}{v^S}
 &=p^S-p^*-\sigma(u^S-u^*)=p^S-p^*+\sigma^2(v^S-v^*),\\
 -\kappa\frac{(\theta^S)'}{\sigma v^S}
 &=c_v(\theta^S-\theta^*)-\frac{1}{2}(u^S-u^*)^2
   -\frac{p^*}{\sigma}(u^S-u^*)\\
 &=c_v(\theta^S-\theta^*)+p^*(v^S-v^*)
      -\frac{1}{2}\sigma^2(v^S-v^*)^2 ,
\end{align*}
where $p^S=\frac{R\theta^S}{v^S}$. The resulting expressions for $(v^S)'$
and $(\theta^S)'$ are smooth functions of $(v^S,\theta^S)$ with
bounded first and second derivatives on the range of the weak shock.
Differentiating twice and using $(u^S)'''=-\sigma(v^S)'''$,
we obtain
\[
 |\partial_\xi^3 q^S|
 \le C\big(|(v^S)''|+|(\theta^S)''|+|(v^S)'|^2+|(\theta^S)'|^2\big)
 \le C\big(\delta_S+(v^S)'\big)(v^S)'\le C(v^S)'.
\]
In particular,
\begin{equation}\label{eq:shock}
 \sum_{j=1}^3|\partial_\xi^j q^S(\xi)|\le C(v^S)'(\xi),\qquad
 0<(v^S)'(\xi)\le C\delta_S^2 e^{-c\delta_S|\xi|},\qquad
 \int_\R(v^S)'\dd\xi=v_+-v^*=\delta_S.
\end{equation}

\subsection{Composite wave and residuals}

Write $\bar p=p(\bar v,\bar\theta)$. In expressions in $x$ a shock
derivative is evaluated at $x-\sigma t-X(t)$
unless its argument is indicated. The profile satisfies
\begin{equation}\label{eq:residual}
 \begin{cases}
 \bar v_t-\bar u_x=-\dot X(v^S)',\\
 \bar u_t+\bar p_x-\mu(\frac{\bar u_x}{\bar v})_x=-\dot X(u^S)'+Q_1,\\
 c_v\bar\theta_t+\bar p\bar u_x-\kappa(\frac{\bar\theta_x}{\bar v})_x
    -\frac{\mu\bar u_x^2}{\bar v}=-c_v\dot X(\theta^S)'+Q_2 .
 \end{cases}
\end{equation}
Put $Q=(Q_1,Q_2)$ and use $Q=Q^R+Q^C+Q^I$ as in
\cite[Section~3]{KVW}, where
\[
 \begin{aligned}
 Q^R&=\left(-\mu(\frac{u^R_x}{v^R})_x,\,
          -\kappa(\frac{\theta^R_x}{v^R})_x-\frac{\mu(u^R_x)^2}{v^R}\right),\\
 Q^C&=\left(u^C_t-\mu(\frac{u^C_x}{v^C})_x,\,-\frac{\mu(u^C_x)^2}{v^C}\right).
 \end{aligned}
\]
The interaction terms in~\cite[Section~3]{KVW} have the form
\begin{align*}
 Q_1^I
 &=\sum_{i=R,C,S}\left[
     (Dp(\bar q)-Dp(q^i))q_x^i
     -\mu\left(\left(\frac{1}{\bar v}-\frac{1}{v^i}\right)u_x^i\right)_x\right],
     \displaybreak[3]\\
 Q_2^I
 &=\sum_{i=R,C,S}\left[
     (\bar p-p^i)u_x^i
     -\kappa\left(\left(\frac{1}{\bar v}-\frac{1}{v^i}\right)\theta_x^i\right)_x
     -\mu\left(\frac{1}{\bar v}-\frac{1}{v^i}\right)(u_x^i)^2\right]\\
 &\quad-\frac{\mu}{\bar v}\sum_{\substack{i,j=R,C,S\\i\ne j}}u_x^i u_x^j ,
\end{align*}
where $p^i=p(v^i,\theta^i)$ and the shock terms are evaluated at
$x-\sigma t-X(t)$.

\begin{lemma}\label{lem:profiles}
If $|\dot X|\le\frac{\sigma}{2}$, then
\begin{align}
 &\|(q^R+q^C)_x\|_\infty+\|(q^R+q^C)_{xx}\|_\infty
           +\|(q^R+q^C)_t\|_\infty
      \le C(1+t)^{-\frac{1}{2}},\label{eq:slow}\\
 &\sup\nolimits_{x\ge\frac{\sigma t}{4}}\sum_{j=0}^2
   |\partial_x^j(q^R+q^C-q_*-q^*)(t,x)|
      \le C(\delta_R+\delta_C)e^{-ct},\label{eq:sep}\\
 &\int_0^T\|Q(t)\|_{H^1}^2\dd t
      \le C(\delta_R+\delta_C),\qquad T>0.\label{eq:QH1}
\end{align}
Moreover,
\begin{equation}\label{eq:lambda}
 \Lambda(t):=\sum_{j=2}^3
   (\|\partial_x^j q^R\|_\infty+\|\partial_x^j q^C\|_\infty)
       +\|q^R_x\|_\infty^2+\|q^C_x\|_\infty^2
 \le C(1+t)^{-1}.
\end{equation}
\end{lemma}

\begin{proof}
The bounds in \eqref{eq:slow} and \eqref{eq:lambda} follow from
\eqref{eq:rare-infty}, \eqref{eq:contact}, and the time derivative
bounds for the contact wave. On $x\ge\frac{\sigma t}{4}$,
\begin{align*}
 &\sum_{j=0}^2|\partial_x^j(q^R+q^C-q_*-q^*)(t,x)|\\
 &\quad\le C\left(\delta_R e^{-2[x-(1+t)w_*]}
                    +\delta_C e^{-\frac{cx^2}{1+t}}\right)
 \le C(\delta_R+\delta_C)e^{-ct}.
\end{align*}
This proves \eqref{eq:sep}.

Take $0<b<\min\{\frac{\sigma}{4},\frac{|w_*|}{2}\}$.
The shift bound gives $\sigma t+X(t)\ge\frac{\sigma t}{2}$.
The end-state bounds through order three give
\[
 \begin{cases}
 \displaystyle\sum_{j=0}^3
 \bigl(|\partial_x^j(q^C-q_*)|+|\partial_x^j(q^S-q^*)|\bigr)
  \le C(\delta_C+\delta_S)e^{-c(t+|x|)},&x\le-bt,\\[2pt]
 \displaystyle\sum_{j=0}^3
 \bigl(|\partial_x^j(q^R-q_*)|+|\partial_x^j(q^S-q^*)|\bigr)
  \le C(\delta_R+\delta_S)e^{-ct},&|x|\le bt,\\[2pt]
 \displaystyle\sum_{j=0}^3
 \bigl(|\partial_x^j(q^R-q_*)|+|\partial_x^j(q^C-q^*)|\bigr)
  \le C(\delta_R+\delta_C)e^{-c(t+|x|)},&x\ge bt .
 \end{cases}
\]
Here the shock arguments are $x-\sigma t-X(t)$, and for $|x|\ge bt$
the Gaussian is bounded by $Ce^{-c(t+|x|)}$.
Differentiating the interaction terms gives
\begin{align*}
 |Q^I|+|Q_x^I|
 &\le C\sum_{i=R,C,S}
       \left(\sum_{j=0}^{2}|\partial_x^j(\bar q-q^i)|\right)
       \left(\sum_{j=1}^{3}|\partial_x^j q^i|\right)\\
 &\quad+C\sum_{\substack{i,k=R,C,S\\i\ne k}}
   (|u_x^i|+|u_{xx}^i|)(|u_x^k|+|u_{xx}^k|)\\
 &\le C(\delta_R+\delta_C)
       \begin{cases}
       e^{-c(t+|x|)},& |x|\ge bt,\\
       e^{-ct},& |x|\le bt .
       \end{cases}
\end{align*}
For the terms with $i=S$ in the first sum, use
$\bar q-q^S=q^R+q^C-q_*-q^*$, whose derivatives through order
two are bounded by $C(\delta_R+\delta_C)$; on $x\ge bt$ they have
the extra factor $e^{-c(t+|x|)}$ from the preceding bounds.
Consequently,
\begin{align*}
 \|Q^I\|_{H^1}
 &\le C(\delta_R+\delta_C)\left(
  \int_{|x|\le bt} e^{-2ct}\dd x+
  \int_{|x|\ge bt} e^{-2c(t+|x|)}\dd x\right)^{\frac{1}{2}}
 \le C(\delta_R+\delta_C)e^{-ct},\\
 \|Q\|_{H^1}
 &\le C\bigl(\|q^R_{xx}\|_2+\|q^R_{xxx}\|_2+\|q^R_x\|_4^2
      +\delta_C(1+t)^{-\frac{5}{4}}+\|Q^I\|_{H^1}\bigr)\\
 &\le C\left(\frac{\delta_R}{1+\delta_R t}
       +\delta_C(1+t)^{-\frac{5}{4}}+(\delta_R+\delta_C)e^{-ct}\right).
\end{align*}
Using \eqref{eq:rare-higher} and the contact bounds as above yields
\[
 \int_0^T\|Q(t)\|_{H^1}^2\dd t
 \le C\int_0^\infty\left(
  \frac{\delta_R^2}{(1+\delta_R t)^2}
  +\delta_C^2(1+t)^{-\frac{5}{2}}
  +(\delta_R+\delta_C)^2e^{-2ct}\right)\dd t
 \le C(\delta_R+\delta_C).
\]
\end{proof}

\section{A priori estimates}\label{sec:energy}

\subsection{Construction of the shift}

We use the form of the shift equation in~\cite[Section~3.2]{KVW}.
In the original Lagrangian coordinate it reads
\begin{equation}\label{eq:shift}
 \begin{aligned}
 \dot X(t)&=-\frac{M}{\delta_S}\int_\R a(x-\sigma t-X(t))
   \left[\frac{\bar p}{\bar v}(v^S)'\phi+(u^S)'\psi
                  +\frac{c_v}{\bar\theta}(\theta^S)'\zeta\right]\dd x,\\
 X(0)&=0,\qquad (\phi,\psi,\zeta)=Z=q-\bar q .
 \end{aligned}
\end{equation}
Here
\[
 a(\xi)=1+\delta_S^{-\frac{1}{4}}(v^S(\xi)-v^*),\quad
 \sigma^*=\sqrt{\frac{\gamma p^*}{v^*}},\quad
 M=\frac{3\gamma(\gamma+1)p^*}{4(v^*)^2(\sigma^*)^3}
       \left(1+\frac{2\kappa(\gamma-1)^2}{\mu R\gamma}\right).
\]
We take $\lambda=\delta_S^{\frac{3}{4}}$ in the weight of
\cite[Section~4.2]{KVW}. The value of $M$ in the final estimate
of~\cite[Section~4.5]{KVW} gives the cancellation
\[
 -(\sigma^*)^2 M
 +\frac{3\gamma(\gamma+1)p^*}{4(v^*)^2\sigma^*}
       \left(1+\frac{2\kappa(\gamma-1)^2}{\mu R\gamma}\right)=0 .
\]
Existence and uniqueness for \eqref{eq:shift} follow
from~\cite[Lemma~3.1]{KVW}. In the small neighborhood,
\begin{equation}\label{eq:shift-bound}
 |\dot X(t)|\le
 \frac{C}{\delta_S}\|Z(t)\|_\infty\int_\R(v^S)'\dd\xi
 \le C\|Z(t)\|_\infty .
\end{equation}
Thus, under a sufficiently small $H^2$ bound, $|\dot X|\le\frac{\sigma}{2}$
and $\frac{\sigma t}{2}\le \sigma t+X(t)\le\frac{3\sigma t}{2}$. The local construction of
\cite[Proposition~3.1]{KVW} is used for the Cauchy problem.

\subsection{The first-order estimate}

We work on an interval $[0,T]$ on which the perturbation is small
in $H^2$. The bounds below are independent of $T$.

Write $\mathcal G^R,\mathcal G^S$ for the two quantities of
\cite[Proposition~3.2]{KVW}:
\[
 \mathcal G^R=\int_\R v^R_x(|\phi|^2+|\zeta|^2)\dd x,\qquad
 \mathcal G^S=\int_\R(v^S)'(x-\sigma t-X(t))|Z|^2\dd x .
\]
Let
\[
 D_1=\mathcal G^R+\mathcal G^S+|\dot X|^2
       +\|Z_x\|_2^2+\|(\psi_{xx},\zeta_{xx})\|_2^2 .
\]
We include the fixed factor $\delta_S$ of the shift dissipation
in the constants. To keep the shock fixed and make the right-hand
side small with $\delta_R,\delta_C$, we use the inequalities in
the proof of~\cite[Proposition~3.2]{KVW} before the wave strengths
are replaced by a common upper bound. With $c_0$ from
\cite[Lemma~4.5]{KVW}, let
\[
 \mathcal W(t)=(1+t)^{-1}\int_\R e^{-\frac{c_0 x^2}{1+t}}|Z|^2\dd x .
\]
The differential inequality following (4.71) of~\cite{KVW},
together with the residual estimates (4.69)--(4.70) and
Lemma~4.2 there, yields
\begin{align}
 &\sup_{0\le t\le T}\|Z(t)\|_2^2+
 \int_0^T(\mathcal G^R+\mathcal G^S+|\dot X|^2+
                                     \|(\psi_x,\zeta_x)\|_2^2)\dd t\notag\\
 &\quad\le C\Big[\|Z(0)\|_2^2+\delta_R^{\frac{4}{3}}+\delta_C^{\frac{4}{3}}
        +(\delta_R+\delta_C)\sup_{0\le t\le T}\|Z(t)\|_2
        +\delta_R^{\frac{1}{3}}\sup_{0\le t\le T}\|Z(t)\|_2^{\frac{2}{3}}
        +\delta_C\int_0^T\mathcal W\dd t\Big]\notag\\
 &\quad\le \frac{1}{2}\sup_{0\le t\le T}\|Z(t)\|_2^2+
      C\Big[\|Z(0)\|_2^2+(\delta_R+\delta_C)^{\frac{1}{2}}
                                      +\delta_C\int_0^T\mathcal W\dd t\Big].
 \label{eq:zero-order}
\end{align}

In the first-order derivative estimates of~\cite[Section~5]{KVW},
the terms with the contact profile have the bound
\begin{align*}
 \int_\R (|q^C_x|^2+|q^C_{xx}|^2)|Z|^2\dd x
 &\le \frac{C\delta_C^2}{1+t}
      \int_\R e^{-\frac{c_0x^2}{1+t}}|Z|^2\dd x
 = C\delta_C^2\mathcal W(t).
\end{align*}
The terms with rarefaction derivatives are bounded by
$C\mathcal G^R+C(1+t)^{-2}\|Z\|_2^2$ by \eqref{eq:rare-infty}.
For the residual terms in (5.8)--(5.12) of~\cite{KVW},
\eqref{eq:QH1} gives
\[
 \int_0^T\|Q(t)\|_2^2\dd t
 \le \int_0^T\|Q(t)\|_{H^1}^2\dd t
 \le C(\delta_R+\delta_C).
\]
Thus, keeping the contact factor in the estimates leading to
(5.8)--(5.12) of~\cite{KVW}, we have
\begin{align*}
 &\sup_{0\le t\le T}\|Z_x(t)\|_2^2+\int_0^T
       \big(\|\phi_x\|_2^2+\|(\psi_{xx},\zeta_{xx})\|_2^2\big)\dd t\\
 &\quad\le C\Big[\|Z(0)\|_{H^1}^2+\sup_{0\le t\le T}\|Z(t)\|_2^2+
     \int_0^T\big(\mathcal G^R+\mathcal G^S+|\dot X|^2
         +\|(\psi_x,\zeta_x)\|_2^2
         +\delta_C^2\mathcal W+\|Q\|_2^2\big)\dd t\Big].
\end{align*}
The estimate of~\cite[Lemma~4.5]{KVW} gives
\[
 \int_0^T\mathcal W\dd t
 \le C\left[1+\sup_{0\le t\le T}\|Z(t)\|_2^2+
                                      \int_0^T D_1\dd t\right].
\]
Consequently, by \eqref{eq:zero-order} and \eqref{eq:QH1},
\begin{align}
 \sup_{0\le t\le T}\|Z(t)\|_{H^1}^2+\int_0^T D_1\dd t
 &\le C\left[\|Z(0)\|_{H^1}^2+(\delta_R+\delta_C)^{\frac{1}{2}}
            +(\delta_C+\delta_C^2)\int_0^T\mathcal W\dd t\right]\notag\\
 &\le C\left[\|Z(0)\|_{H^1}^2+(\delta_R+\delta_C)^{\frac{1}{2}}\right],
 \label{eq:first-order}
\end{align}
where the contact strength is sufficiently small for the absorption.

\subsection{The second-order estimate}

Subtracting \eqref{eq:residual} from \eqref{eq:nsf} and
\eqref{eq:temperature}, we obtain
\begin{equation}\label{eq:perturb}
 \begin{aligned}
 \phi_t&=\psi_x+\dot X(v^S)',\\
 \psi_t&=-(p-\bar p)_x
     +\mu\left(\frac{u_x}{v}-\frac{\bar u_x}{\bar v}\right)_x
     +\dot X(u^S)'-Q_1,\\
 c_v\zeta_t&=-p\psi_x-(p-\bar p)\bar u_x
     +\kappa\left(\frac{\theta_x}{v}-\frac{\bar\theta_x}{\bar v}\right)_x\\
 &\quad+\mu\left(\frac{u_x^2}{v}-\frac{\bar u_x^2}{\bar v}\right)
        +c_v\dot X(\theta^S)'-Q_2 .
 \end{aligned}
\end{equation}
\begin{proposition}\label{prop:energy}
Under \eqref{eq:small}, the solution satisfies
\begin{equation}\label{eq:energy}
 \sup_{t\ge0}\|Z(t)\|_{H^2}^2+
 \int_0^\infty\big[D_1+\|\phi_{xx}\|_2^2+
                         \|(\psi_{xxx},\zeta_{xxx})\|_2^2\big]\dd t
 \le C_S\big[\|Z(0)\|_{H^2}^2+(\delta_R+\delta_C)^{\frac{1}{2}}\big].
\end{equation}
\end{proposition}

\begin{proof}
We work on $[0,T]$ with $\sup_{0\le t\le T}\|Z(t)\|_{H^2}$
sufficiently small. For example, the second derivative of the
viscous flux is
\begin{align*}
 \left(\frac{u_x}{v}-\frac{\bar u_x}{\bar v}\right)_{xx}
 &=\frac{\psi_{xxx}}{v}-\frac{u_x}{v^2}\phi_{xx}
   +\bar u_{xxx}\left(\frac{1}{v}-\frac{1}{\bar v}\right)
   -\bar v_{xx}\left(\frac{u_x}{v^2}-\frac{\bar u_x}{\bar v^2}\right)\\
 &\quad-2\left(\frac{v_xu_{xx}}{v^2}
       -\frac{\bar v_x\bar u_{xx}}{\bar v^2}\right)
       +2\left(\frac{u_x v_x^2}{v^3}
       -\frac{\bar u_x\bar v_x^2}{\bar v^3}\right).
\end{align*}
The same product rule for the heat flux and the pressure gives
\begin{align}
 &\|\left(\frac{u_x}{v}-\frac{\bar u_x}{\bar v}\right)_{xx}-\frac{\psi_{xxx}}{v}+\frac{u_x\phi_{xx}}{v^2}\|_2
   +\|(p-\bar p)_{xx}+\frac{R\theta}{v^2}\phi_{xx}-\frac{R}{v}\zeta_{xx}\|_2\notag\\
 &\quad+\|\left(\frac{\theta_x}{v}-\frac{\bar\theta_x}{\bar v}\right)_{xx}-\frac{\zeta_{xxx}}{v}+\frac{\theta_x\phi_{xx}}{v^2}\|_2\notag\\
 &\quad\le C\big[\|Z_x\|_2+\|(\psi_{xx},\zeta_{xx})\|_2+
               \|(|\bar q_{xx}|+|\bar q_{xxx}|+|\bar q_x|^2)|Z|\|_2\big]\notag\\
 &\quad\le C\big[D_1^{\frac{1}{2}}+\Lambda(t)\|Z\|_2\big].
 \label{eq:flux-two}
\end{align}
The factors with shock derivatives in \eqref{eq:flux-two} are
bounded in $L^2$ by $C(\mathcal G^S)^{\frac{1}{2}}$ by \eqref{eq:shock}.
For the other profile derivatives we use \eqref{eq:lambda}. For the
pressure work and the viscous heating we obtain
\begin{align}
 &\|(p\psi_x+(p-\bar p)\bar u_x)_{xx}-p\psi_{xxx}\|_2+
       \|\left(\frac{u_x^2}{v}-\frac{\bar u_x^2}{\bar v}\right)_{xx}-\frac{2u_x\psi_{xxx}}{v}+\frac{u_x^2\phi_{xx}}{v^2}\|_2\notag\\
 &\quad\le C\big[D_1^{\frac{1}{2}}+\Lambda\|Z\|_2+
   (\|\psi_x\|_\infty+\|\bar u_x\|_\infty)\|\phi_{xx}\|_2+
                                 \|\psi_{xx}\|_4^2\big].
 \label{eq:heat}
\end{align}
Here the differentiated pressure terms are
$2p_x\psi_{xx}+p_{xx}\psi_x+(p-\bar p)_{xx}\bar u_x+
2(p-\bar p)_x\bar u_{xx}+(p-\bar p)\bar u_{xxx}$.

Differentiate \eqref{eq:perturb} twice in $x$ and multiply its
equations by $\frac{R\theta}{v^2}\phi_{xx}$, $\psi_{xx}$ and $\frac{\zeta_{xx}}{\theta}$,
respectively. Integration by parts in the momentum equation
gives the cancellation
\[
 \frac{R\theta}{v^2}\phi_{xx}\psi_{xxx}+(-\frac{R\theta}{v^2}\phi_{xx}+\frac{R}{v}\zeta_{xx})\psi_{xxx}
                       -\frac{p}{\theta}\zeta_{xx}\psi_{xxx}=0 .
\]
The diffusion terms give
$\int_\R[\frac{\mu\psi_{xxx}^2}{v}+\frac{\kappa\zeta_{xxx}^2}{v\theta}]\dd x$.
The other density terms in the diffusion fluxes have the factors
$u_x$ and $\theta_x$ in \eqref{eq:flux-two}.
Using $v_t=u_x$ and \eqref{eq:temperature}, we have
\[
 |\left(\frac{R\theta}{v^2}\right)_t|+|\left(\frac{1}{\theta}\right)_t|
 \le C(|u_x|+|\theta_xv_x|+|u_x|^2+|\bar\theta_{xx}|+|\zeta_{xx}|).
\]
The remaining products involving second derivatives satisfy
\begin{align*}
 &\int_\R\big(|\zeta_{xx}|\phi_{xx}^2+|\zeta_{xx}|^3+
                       |\psi_{xx}|^2|\zeta_{xx}|\big)\dd x\\
 &\quad\le C\|\zeta_{xx}\|_2^{\frac{1}{2}}\|\zeta_{xxx}\|_2^{\frac{1}{2}}
                   (\|\phi_{xx}\|_2^2+\|\zeta_{xx}\|_2^2)
          +C\|\zeta_{xx}\|_2\|\psi_{xx}\|_2^{\frac{3}{2}}\|\psi_{xxx}\|_2^{\frac{1}{2}}\\
 &\quad\le\eta\|(\psi_{xxx},\zeta_{xxx})\|_2^2+
       C_\eta\|Z\|_{H^2}^{\frac{4}{3}}\big(\|\phi_{xx}\|_2^2+
                                      \|(\psi_{xx},\zeta_{xx})\|_2^2\big).
\end{align*}
Integrating the $Q$ terms by parts and applying
\eqref{eq:flux-two}--\eqref{eq:heat} gives
\begin{align}
 &\frac{\mathrm d}{\mathrm dt}\int_\R
     \left(\frac{R\theta\phi_{xx}^2}{2v^2}
          +\frac{\psi_{xx}^2}{2}
          +\frac{c_v\zeta_{xx}^2}{2\theta}\right)\dd x
         +c\|(\psi_{xxx},\zeta_{xxx})\|_2^2\notag\\
 &\quad\le C\big(D_1^{\frac{1}{2}}+\Lambda\|Z\|_2+\|Q_x\|_2
           +(\|q_x\|_\infty+\|\bar q_x\|_\infty)\|\phi_{xx}\|_2\big)
      \big(\|(\psi_{xxx},\zeta_{xxx})\|_2+
            \|(\psi_{xx},\zeta_{xx})\|_2\big)\notag\displaybreak[3]\\
 &\qquad+C|\dot X|\,\|\phi_{xx}\|_2
       +C\big(\|q_x\|_\infty+\|\bar\theta_{xx}\|_\infty\big)
                      \|Z_{xx}\|_2^2
       +C\int_\R\big(|\zeta_{xx}|\phi_{xx}^2+|\zeta_{xx}|^3+
                             |\psi_{xx}|^2|\zeta_{xx}|\big)\dd x\notag\\
 &\quad\le \frac{c}{2}\|(\psi_{xxx},\zeta_{xxx})\|_2^2
       +\eta\|\phi_{xx}\|_2^2
       +C_\eta(D_1+\Lambda^2\|Z\|_2^2+\|Q_x\|_2^2)
 \label{eq:second-energy}
\end{align}
for any fixed sufficiently small $\eta>0$, with the wave strengths
and the $H^2$ bound taken sufficiently small. The term
$C|\dot X|\|\phi_{xx}\|_2$ comes from the differentiated mass equation.

The first two equations in
\eqref{eq:perturb} give
\begin{align*}
 \frac{\mathrm d}{\mathrm dt}\int_\R\psi_x\phi_{xx}\dd x
 &=\int_\R\left[-(p-\bar p)_{xx}
            +\mu\left(\frac{u_x}{v}-\frac{\bar u_x}{\bar v}\right)_{xx}
       +\dot X(u^S)''-Q_{1,x}\right]\phi_{xx}\dd x\\
 &\quad-\|\psi_{xx}\|_2^2+\dot X\int_\R\psi_x(v^S)'''\dd x .
\end{align*}
It follows from \eqref{eq:flux-two} that
\begin{align}
 &-\frac{\mathrm d}{\mathrm dt}\int_\R\psi_x\phi_{xx}\dd x
      +\int_\R \frac{R\theta}{v^2}\phi_{xx}^2\dd x\notag\\
 &\quad\le C\|\phi_{xx}\|_2
       \big[\|\psi_{xxx}\|_2+D_1^{\frac{1}{2}}+\Lambda\|Z\|_2+\|Q_x\|_2\big]
          +C\|u_x\|_\infty\|\phi_{xx}\|_2^2+C D_1\notag\\
 &\le\frac{1}{2}\int_\R \frac{R\theta}{v^2}\phi_{xx}^2\dd x+
       C[\|\psi_{xxx}\|_2^2+D_1+\Lambda^2\|Z\|_2^2+\|Q_x\|_2^2].
 \label{eq:density-test}
\end{align}
Absorb the third-derivative term on the right of \eqref{eq:second-energy}
and the density term on the right of \eqref{eq:density-test}. Add a
sufficiently small multiple $\rho$ of the latter estimate to
the former, and then take $\eta$ sufficiently small.
We obtain
\[
 \begin{aligned}
 &\frac{\mathrm d}{\mathrm dt}\int_\R
     \left(\frac{R\theta\phi_{xx}^2}{2v^2}
          +\frac{\psi_{xx}^2}{2}
          +\frac{c_v\zeta_{xx}^2}{2\theta}
          -\rho\psi_x\phi_{xx}\right)\dd x
     +c\bigl(\|\phi_{xx}\|_2^2+\|(\psi_{xxx},\zeta_{xxx})\|_2^2\bigr)\\
 &\qquad\le C\bigl(D_1+\Lambda^2\|Z\|_2^2+\|Q_x\|_2^2\bigr).
 \end{aligned}
\]
By Young's inequality,
$
\rho\left|\int_\R\psi_x\phi_{xx}\dd x\right|
\le\int_\R\frac{R\theta\phi_{xx}^2}{4v^2}\dd x
      +C\rho^2\|\psi_x\|_2^2.
$
Integrating and using \eqref{eq:first-order} and
Lemma~\ref{lem:profiles}, we have
\begin{align*}
 &\sup_{0\le t\le T}\|Z(t)\|_{H^2}^2+
     \int_0^T[D_1+\|\phi_{xx}\|_2^2+\|(\psi_{xxx},\zeta_{xxx})\|_2^2]\dd t\\
 &\quad\le C\left[\|Z(0)\|_{H^2}^2+\sup_{0\le t\le T}\|Z(t)\|_{H^1}^2+
         \int_0^T(D_1+\|Q_x\|_2^2)\dd t+
                   \|\Lambda\|_{L^2(0,T)}^2\sup_{0\le t\le T}\|Z(t)\|_2^2\right]\\
 &\quad\le C_S\big[\|Z(0)\|_{H^2}^2+(\delta_R+\delta_C)^{\frac{1}{2}}\big].
\end{align*}

The estimate applies first to smooth data. It keeps the $H^2$
norm inside the small neighborhood, so positivity and the
continuation criterion of~\cite[Section~3]{KVW} give a global
solution. For $H^2$ data, the estimates pass to the limit by
smooth approximation, local compactness and weak lower
semicontinuity. On each finite interval, \eqref{eq:perturb} and
the bounds give
\[
 \phi_t\in L^2(0,T;H^2),\qquad
 (\psi_t,\zeta_t)\in L^2(0,T;H^1),\qquad
 (\psi,\zeta)\in L^2(0,T;H^3).
\]
Continuity of $\phi$ into $H^2$ follows by integration in time.
The last two inclusions give $(\psi,\zeta)\in C([0,T];H^2)$
by the energy identity with $H^2$ as the pivot space. The
right-hand side of \eqref{eq:shift} is continuous by these
bounds and the continuity of $X$, hence $X\in C^1$.
Uniqueness follows from~\cite[Proposition~3.1]{KVW}. Letting
$T\to\infty$ gives \eqref{eq:energy}.
\end{proof}

\section{Comparison estimates}\label{sec:comparison}

\subsection{An exponentially weighted estimate}

\begin{lemma}\label{lem:local}
Let $q,\widetilde q$ be strong solutions of \eqref{eq:nsf} on
$[s,s+\ell]$, $\ell\ge1$. Suppose that the states lie in a fixed
compact positive range, their first spatial derivatives are uniformly
bounded, and
\begin{equation}\label{eq:local-input}
 \int_s^{s+\ell}\|\theta_{xx}(\tau)\|_\infty\dd\tau\le C\ell .
\end{equation}
If the solutions agree at $s$ on $|x-x_0|\le R_0$, there is a
constant $C_0>0$, depending only on these bounds, such that
\begin{equation}\label{eq:local}
 \sup_{\substack{s\le\tau\le s+\ell\\|x-x_0|\le r}}
       |q(\tau,x)-\widetilde q(\tau,x)|
 \le C\exp\bigl(-\frac{R_0-r-C_0\ell}{4}\bigr)
\end{equation}
whenever $R_0-r-C_0\ell\ge1$.
\end{lemma}

\begin{proof}
Let $(f,g,h)=q-\widetilde q$ and
$\omega(x)=\exp(-\sqrt{1+(x-x_0)^2})$. We use the weighted energy
\[
 \mathcal E=\frac{1}{2}\int_\R\omega
       \left(\frac{R\theta}{v^2}f^2+g^2+\frac{c_vh^2}{\theta}\right)\dd x.
\]
The pressure and diffusion terms satisfy
\begin{align*}
 (p-\widetilde p)_x&=-\frac{R\theta}{v^2}f_x+\frac{R}{v}h_x
                      +(Dp(q)-Dp(\widetilde q))\widetilde q_x,\\
 \frac{u_x}{v}-\frac{\widetilde u_x}{\widetilde v}
 &=\frac{g_x}{v}-\frac{\widetilde u_x f}{v\widetilde v},
 &\frac{\theta_x}{v}-\frac{\widetilde\theta_x}{\widetilde v}
 &=\frac{h_x}{v}-\frac{\widetilde\theta_x f}{v\widetilde v}.
\end{align*}
Multiply the difference equations by
$\frac{\omega R\theta}{v^2}f$, $\omega g$, and
$\frac{\omega h}{\theta}$. The first-order terms and the time
derivatives of the coefficients are
\begin{align*}
 \omega\left(\frac{R\theta}{v^2}(fg_x+gf_x)
       -\frac{R}{v}(gh_x+hg_x)\right)
 &=\omega\left(\frac{R\theta}{v^2}(fg)_x-\frac{R}{v}(gh)_x\right),\\
 \left|\left(\frac{R\theta}{v^2}\right)_t\right|
 +\left|\left(\frac{1}{\theta}\right)_t\right|
 &\le C(1+|\theta_{xx}|),
\end{align*}
by $v_t=u_x$ and \eqref{eq:temperature}. Using the identities
above and $|\frac{u_x^2}{v}-\frac{\widetilde u_x^2}{\widetilde v}|\le C(|g_x|+|f|)$,
we obtain, by integration by parts and Young's inequality,
\begin{align*}
 &\mathcal E'+\int_\R\omega\left(
                  \frac{\mu g_x^2}{v}+\frac{\kappa h_x^2}{v\theta}
                          \right)\dd x\\
 &\quad\le \int_\R\big[\left(\frac{\omega R}{v}\right)_x gh-\left(\frac{\omega R\theta}{v^2}\right)_x fg\big]\dd x\\
 &\qquad+C\int_\R\omega\big[(1+|\theta_{xx}|)|(f,g,h)|^2
                        +|(f,g,h)|\,|(g_x,h_x)|\big]\dd x
       \displaybreak[3]\\
 &\quad\le C(1+\|\theta_{xx}\|_\infty)\mathcal E+
       \frac{1}{2}\int_\R\omega\left(
                  \frac{\mu g_x^2}{v}+\frac{\kappa h_x^2}{v\theta}
                          \right)\dd x,
\end{align*}
where $|\omega_x|+|\left(\frac{\omega R\theta}{v^2}\right)_x|+|\left(\frac{\omega R}{v}\right)_x|\le C\omega$.
The integrals on $\R$ are justified by spatial cutoffs.

Since $(f,g,h)(s,x)=0$ for $|x-x_0|\le R_0$, Gronwall's
inequality and \eqref{eq:local-input} give, for $s\le\tau\le s+\ell$,
\begin{align*}
 \mathcal E(\tau)
 &\le\mathcal E(s)\exp\!\left(
         C\int_s^\tau (1+\|\theta_{xx}(t)\|_\infty)\dd t\right)\\
 &\le C e^{C_0\ell}\int_{|x-x_0|>R_0}e^{-|x-x_0|}\dd x
 \le C e^{-R_0+C_0\ell}.
\end{align*}
For $|x-x_0|\le r$ and $|y-x|\le1$, we have
$\omega(y)\ge e^{-r-2}$. By averaging over $[x-1,x+1]$ and
using the uniform first-derivative bound, we get
\begin{align*}
 &|(f,g,h)(\tau,x)|^2\\
 &\quad\le \frac{1}{2}\int_{x-1}^{x+1}|(f,g,h)(\tau,y)|^2\dd y
   +2\int_{x-1}^{x+1}|(f,g,h)(\tau,y)|\,|(f_y,g_y,h_y)(\tau,y)|\dd y\\
 &\quad\le C e^r\mathcal E(\tau)+C(e^r\mathcal E(\tau))^{\frac{1}{2}}
 \le C e^{-\frac{R_0-r-C_0\ell}{2}}.
\end{align*}
Taking square roots proves \eqref{eq:local}.
\end{proof}

For the composite solution, the profile bounds and \eqref{eq:energy}
give
\begin{align*}
 \int_s^{s+\ell}\|\theta_{xx}\|_\infty\dd\tau
 &\le\int_s^{s+\ell}\|\bar\theta_{xx}\|_\infty\dd\tau+
  C\ell^{\frac{1}{2}}\left(\int_s^{s+\ell}\|\zeta_{xx}\|_{H^1}^2
                      \dd\tau\right)^{\frac{1}{2}}\\
 &\le C\ell+C\ell^{\frac{1}{2}}\le C\ell,
\end{align*}
which is \eqref{eq:local-input}.

\subsection{Single-shock and constant-state estimates}

For the fixed shock $q^S$, use the shift equation
\[
 Y(0)=0,\qquad
 \dot Y(t)=\mathcal L\big(\widetilde q(t,\sigma t+Y(t)+\cdot)-q^S\big),
\]
where
\begin{equation}\label{eq:shock-functional}
 \mathcal L(f,g,h)=-\frac{M}{\delta_S}\int_\R a
 \left(\frac{p^S}{v^S}(v^S)'f+(u^S)'g+
                      \frac{c_v}{\theta^S}(\theta^S)'h\right)\dd\xi .
\end{equation}
The weight $a$ and coefficient $M$ are those in \eqref{eq:shift}.
For the single shock, they agree with the choice in
\cite[(2.7)--(2.9)]{single} with left state $q^*$.
For the datum $q^S+f$ with $\|f\|_{H^2}$ sufficiently small,
\cite[Theorem~2.1 and (2.12)]{single} gives
\begin{equation}\label{eq:single-input}
 \sup_{t\ge0}(1+t)^{\frac{1}{4}}
 \|\widetilde q(t,\cdot)-q^S(\cdot-\sigma t-Y(t))\|_\infty
 \le C_S\|f\|_{H^2}.
\end{equation}
The cited estimate also includes $|\dot Y(t)|$ on the left, with
$C_S$ independent of $t$ and $f$.

\begin{lemma}\label{lem:constant}
Let
$q_c=(v_c,u_c,\theta_c)$ range in a fixed compact set of positive states. For $\|f\|_{H^2}$ sufficiently small,
uniformly in $q_c$, the solution with initial datum
$q_c+f$ satisfies
\begin{equation}\label{eq:constant}
 \sup_{t\ge0}\Big(\|q(t)-q_c\|_{H^2}
       +(1+t)^{\frac{1}{4}}\|q(t)-q_c\|_\infty\Big)\le C\|f\|_{H^2}.
\end{equation}
\end{lemma}

\begin{proof}
First take $q_c=q_+$. Let $q_L$ be the solution with
datum $q^S(x)+f(x-L)$ and let $Y_L$ be its shift in
\eqref{eq:single-input}. For each finite $T$, the profiles
$q^S(\cdot+L-\sigma t-Y_L(t))$ converge to $q_+$ with their
spatial derivatives on compact intervals, uniformly for $0\le t\le T$.
Indeed, the single-shock estimate bounds $\dot Y_L$ independently
of $L$. The $H^2$ bound of~\cite[Proposition~3.2 and (3.6)]{single}
and the equations \eqref{eq:nsf} give
local bounds in $L^\infty(0,T;L^2)$ for the time derivatives
of $q_L(t,\cdot+L)$. Thus a subsequence converges in
$C([0,T];H^1_{\rm loc})$ to the solution $q$ with datum $q_++f$.
The fluxes pass to the limit in the distributional equations by
this convergence and the positive lower bound for $v$.
Along this subsequence, for each compact interval $I$,
\begin{align*}
 &\|q(t)-q_+\|_{H^2(I)}+(1+t)^{\frac{1}{4}}\|q(t)-q_+\|_{L^\infty(I)}\\
 &\quad\le\liminf_{L\to\infty}\bigg(
     \|q_L(t)-q^S(\cdot-\sigma t-Y_L(t))\|_{H^2}
   +(1+t)^{\frac{1}{4}}
       \|q_L(t)-q^S(\cdot-\sigma t-Y_L(t))\|_\infty\bigg)\\
 &\quad\le C_S\|f\|_{H^2}.
\end{align*}
The derivative dissipation in~\cite[(3.6)]{single} also passes to the
limit by weak lower semicontinuity on compact intervals. Uniqueness
in the small perturbation class identifies the limits on different
time intervals. Letting $I$ increase to $\R$ gives
\eqref{eq:constant} for $q_c=q_+$.

For general $q_c$, use the rescaling
\begin{align*}
 y&=\sqrt{\frac{\theta_c}{\theta_+}}\,x,\qquad
 \tau=\frac{\theta_c v_+}{\theta_+v_c}\,t,\\
 (\widehat v,\widehat u,\widehat\theta)
 &=\left(\frac{v_+}{v_c}v,\ 
         u_++\sqrt{\frac{\theta_+}{\theta_c}}(u-u_c),\
         \frac{\theta_+}{\theta_c}\theta\right).
\end{align*}
Substitution preserves \eqref{eq:nsf} with the same coefficients
and maps $q_c$ to $q_+$. The time and norm equivalences are
uniform on the given compact set.
\end{proof}

\section{Proof of Theorem~\ref{thm:decay}}\label{sec:proof}

Proposition~\ref{prop:energy} gives the global solution. By
\eqref{eq:perturb}, \eqref{eq:energy}, and the residual estimates in
the proof of Lemma~\ref{lem:profiles},
\begin{align*}
 \|Z_t(t)\|_2
 &\le C\big(\|Z(t)\|_{H^2}+|\dot X(t)|+\|Q(t)\|_2\big)\le C_S,\\
 \|Z(t)-Z(s)\|_\infty^2
 &\le C\|Z(t)-Z(s)\|_2\|Z_x(t)-Z_x(s)\|_2\le C_S|t-s|.
\end{align*}
By \eqref{eq:shift-bound} and the Sobolev inequality, for $0\le s<T$,
\begin{align*}
 &\int_s^T\big(\|Z(\tau)\|_\infty+|\dot X(\tau)|\big)^4\dd\tau\\
 &\le C\sup_{\tau\ge0}\|Z(\tau)\|_2^2
             \int_s^T\|Z_x(\tau)\|_2^2\dd\tau
 \le C_S\big(\|Z(0)\|_{H^2}^2+(\delta_R+\delta_C)^{\frac{1}{2}}\big)^2 .
\end{align*}
Taking $s=0$ and $T\to\infty$ gives $\|Z\|_\infty\in L^4(0,\infty)$.
Uniform continuity and \eqref{eq:shift-bound} give
$\|Z(t)\|_\infty+|\dot X(t)|\to0$. Taking $s=t$, $T=2t$ and dividing
by $t$ proves the time-averaged bound since $\frac{1}{t}\le\frac{2}{1+t}$
for $t\ge1$.
Fix a smooth cutoff $0\le\chi\le1$ equal to one on $[-1,1]$
and zero outside $[-2,2]$. The estimate \eqref{eq:energy}
allows us to choose $\eps_S$ below the smallness thresholds of
the comparison estimates. We first take $t$ large enough that
the comparison intervals have length at least one.

\subsection{Estimate near the shock}

We compare on $[s,t]$, where $\ell=\nu t$ and $s=t-\ell$.
For $2A_0\nu\le\frac{\sigma(1-\nu)}{4}$, the support of the cutoff
below lies in $x\ge\frac{\sigma s}{4}$, since
\[
 \sigma s+X(s)-2A_0\ell
 \ge\frac{\sigma s}{2}-2A_0\ell
 \ge\frac{\sigma s}{4}.
\]
Let $\widetilde q$ solve \eqref{eq:nsf} with the datum at time $s$
\[
 \widetilde q(s,x)=q^S(x-\sigma s-X(s))
 +\chi\!\left(\frac{x-\sigma s-X(s)}{A_0\ell}\right)
             \big(q(s,x)-q^S(x-\sigma s-X(s))\big).
\]
By \eqref{eq:sep} and \eqref{eq:energy},
\begin{align*}
 &\|\widetilde q(s)-q^S(\cdot-\sigma s-X(s))\|_{H^2}\\
 &\quad\le C\|Z(s)\|_{H^2}+C\left(
    \int_{|x-\sigma s-X(s)|\le2A_0\ell}
     \sum_{j=0}^2|\partial_x^j(q^R+q^C-q_*-q^*)(s,x)|^2\dd x
     \right)^{\frac{1}{2}}\\
 &\quad\le C\|Z(s)\|_{H^2}+C(A_0\ell)^{\frac{1}{2}}e^{-cs}
 \le C_S\eps_S+C(A_0\ell)^{\frac{1}{2}}e^{-cs}.
\end{align*}
For small $\eps_S$ and large $t$, apply \eqref{eq:single-input}
with initial time $s$ and spatial origin $\sigma s+X(s)$.
Write the comparison shift as $Y$,
with $Y(s)=X(s)$ and
$\dot Y(\tau)=\mathcal L(\widetilde q(\tau,\sigma\tau+Y(\tau)+\cdot)-q^S)$.
Since $\sigma\tau+X(\tau)+\xi\ge\frac{\sigma\tau}{4}$ for
$|\xi|\le\frac{\sigma\tau}{4}$,
\begin{align}
 &\left|\dot X(\tau)
       -\mathcal L(q(\tau,\sigma\tau+X(\tau)+\cdot)-q^S)\right|\notag\\
 &\quad\le Ce^{-c\tau}\int_{|\xi|\le\frac{\sigma\tau}{4}}(v^S)'(\xi)\dd\xi
       +C\int_{|\xi|>\frac{\sigma\tau}{4}}(v^S)'(\xi)\dd\xi
 \le Ce^{-c\tau}.
 \label{eq:shift-error}
\end{align}
The solutions coincide at $s$ for
$|x-\sigma s-X(s)|\le A_0\ell$. By \eqref{eq:shift-bound},
\[
 |\sigma(\tau-s)+X(\tau)-X(s)|
 \le\int_s^\tau(\sigma+|\dot X(r)|)\dd r\le2\sigma\ell.
\]
Use Lemma~\ref{lem:local} with center $\sigma s+X(s)$,
$R_0=A_0\ell$ and $r=(\frac{A_0}{2}+2\sigma)\ell$.
The bounds for $q$ and $\widetilde q$ in Lemma~\ref{lem:local},
including \eqref{eq:local-input}, follow from \eqref{eq:energy}
and~\cite[(3.6)]{single}.
Subtracting the shift equations and using the uniform bound on
$\widetilde q_x$ gives
\begin{align}
 |\dot X(\tau)-\dot Y(\tau)|
 &\le C|X(\tau)-Y(\tau)|+Ce^{-cs}\notag\\
 &\quad+C\int_{|\xi|\le\frac{A_0\ell}{2}}(v^S)'(\xi)
   |(q-\widetilde q)(\tau,\sigma\tau+X(\tau)+\xi)|\dd\xi
       +C\int_{|\xi|>\frac{A_0\ell}{2}}(v^S)'(\xi)\dd\xi\notag\\
 &\le C|X(\tau)-Y(\tau)|+
  C\left(e^{-cs}+e^{-\frac{(\frac{A_0}{2}-2\sigma-C_0)\ell}{4}}
                   +e^{-\frac{cA_0\ell}{2}}\right).
 \label{eq:D}
\end{align}
In the first inequality we used \eqref{eq:shift-error}. Choose $A_0$
large and then $\nu$ small so that the three exponentials are at
most $e^{-(C+3)\ell}$, where $C$ is the coefficient of
$|X-Y|$ in \eqref{eq:D}. Gronwall's inequality yields
\[
 \sup_{s\le\tau\le t}|X(\tau)-Y(\tau)|
 \le C\int_s^t e^{C(t-\tau)}e^{-(C+3)\ell}\dd\tau
 \le C\ell e^{-3\ell}.
\]

For $|x-\sigma t-X(t)|\le\ell$ we have
$|x-\sigma s-X(s)|\le(1+2\sigma)\ell$ and $x\ge\frac{\sigma t}{4}$.
Using \eqref{eq:local}, \eqref{eq:single-input} and \eqref{eq:sep},
\begin{align}
 &|q(t,x)-\bar q(t,x)|\notag\\
 &\quad\le |q(t,x)-\widetilde q(t,x)|
            +|\widetilde q(t,x)-q^S(x-\sigma t-Y(t))|\notag\\
 &\qquad\quad+|q^S(x-\sigma t-Y(t))-q^S(x-\sigma t-X(t))|
              +|(q^R+q^C-q_*-q^*)(t,x)|\notag\\
 &\quad\le C\big(e^{-c\ell}+(1+\ell)^{-\frac{1}{4}}
                      +\ell e^{-3\ell}+e^{-ct}\big)
 \le C(1+t)^{-\frac{1}{4}} .
 \label{eq:near-final}
\end{align}
The shift bound follows from the same estimate,
\begin{align*}
 |\dot X(t)|
 &\le C\int_\R(v^S)'(\xi)|Z(t,\sigma t+X(t)+\xi)|\dd\xi\\
 &\le C\sup_{|x-\sigma t-X(t)|\le\ell}|Z(t,x)|
                       +C\int_{|\xi|>\ell}(v^S)'(\xi)\dd\xi
 \le C\big((1+t)^{-\frac{1}{4}}+e^{-c\ell}\big)
 \le C(1+t)^{-\frac{1}{4}} .
\end{align*}
This proves \eqref{eq:shock-rate} for large $t$.

\subsection{Estimates away from the shock}

For $|x_0-\sigma t-X(t)|>\nu t$, compare the solutions on
$[t-\ell,t]$ with the same $A_0>C_0+4$ as in Section~5.1. If
$1\le\ell\le\min\{\frac{t}{2},\frac{\nu t}{4A_0+3\sigma}\}$, then for
$t-\ell\le\tau\le t$ and $|x-x_0|\le2A_0\ell$,
\begin{align*}
 |x-\sigma\tau-X(\tau)|
 &\ge |x_0-\sigma t-X(t)|-|x-x_0|
                   -\int_\tau^t(\sigma+|\dot X(r)|)\dd r\\
 &\ge \nu t-\left(2A_0+\frac{3\sigma}{2}\right)\ell
 \ge \frac{\nu t}{2}
\end{align*}
by \eqref{eq:shift-bound}. If $|x_0|\ge\eta t$ and
$\ell\le\frac{\eta t}{4A_0}$, then
$|x|\ge|x_0|-2A_0\ell\ge\frac{\eta t}{2}$ on the same interval.
Since $1+\tau\ge\frac{1+t}{2}$, the bounds of Section~2 give
\begin{align*}
 &|\bar q_x(\tau,x)|+|\bar q_{xx}(\tau,x)|+|\bar q_t(\tau,x)|\\
 &\quad\le C\left((1+t)^{-1}
           +(1+t)^{-\frac{1}{2}}e^{-\frac{c x^2}{1+t}}
           +e^{-ct}\right)
 \le\begin{cases}
 C(1+t)^{-\frac{1}{2}},& |x_0-\sigma t-X(t)|>\nu t,\\
 C_\eta(1+t)^{-1},& |x_0-\sigma t-X(t)|>\nu t,\quad |x_0|\ge\eta t .
 \end{cases}
\end{align*}

At $s=t-\ell$ use the perturbation
\[
 f(x)=\chi\!\left(\frac{x-x_0}{A_0\ell}\right)
                 \big(q(s,x)-\bar q(s,x_0)\big).
\]
The product rule and the identity
$\bar q(s,x)-\bar q(s,x_0)=\int_{x_0}^{x}\bar q_x(s,y)\dd y$
give, for $\ell\ge1$,
\begin{align}
 \|f\|_{H^2}
 &\le C\|Z(s)\|_{H^2}+C\ell^{\frac{1}{2}}
    \sup_{|x-x_0|\le2A_0\ell}
    \left(\left|\int_{x_0}^{x}\bar q_x(s,y)\dd y\right|
                   +|\bar q_x(s,x)|+|\bar q_{xx}(s,x)|\right)\notag\\
 &\le C\|Z(s)\|_{H^2}+C\ell^{\frac{3}{2}}
     \sup_{|x-x_0|\le2A_0\ell} (|\bar q_x(s,x)|+|\bar q_{xx}(s,x)|)\notag\\
 &\le C\|Z(s)\|_{H^2}
    +\begin{cases}
       C\ell^{\frac{3}{2}}(1+t)^{-\frac{1}{2}},
               &\text{for the bound on the whole line},\\
       C_\eta\ell^{\frac{3}{2}}(1+t)^{-1},
               &\text{for the bound on } |x_0|\ge\eta t .
      \end{cases}
 \label{eq:constant-cutoff}
\end{align}
Take $\ell=\eps_0(1+t)^{\frac{1}{3}}$ and
$\ell=\eps_0(1+t)^{\frac{2}{3}}$ in the two cases, respectively.
For large $t$ both lengths satisfy the restrictions above. By
\eqref{eq:energy} and \eqref{eq:constant-cutoff}, choosing $\eps_0$
and $\eps_S$ sufficiently small gives the solution $\widetilde q$ of
Lemma~\ref{lem:constant} with datum $\bar q(s,x_0)+f$, which equals
$q(s,x)$ on $|x-x_0|\le A_0\ell$.

Taking $R_0=A_0\ell$ and $r=1$ in \eqref{eq:local} and using
\eqref{eq:constant}, we obtain
\begin{align}
 &|q(t,x_0)-\bar q(t,x_0)|\notag\\
 &\quad\le |q(t,x_0)-\widetilde q(t,x_0)|
      +|\widetilde q(t,x_0)-\bar q(s,x_0)|
      +\int_s^t|\bar q_t(\tau,x_0)|\dd\tau\notag\\
 &\quad\le C\big(e^{-c\ell}+(1+\ell)^{-\frac{1}{4}}\|f\|_{H^2}\big)
             +\ell\sup_{s\le\tau\le t}|\bar q_t(\tau,x_0)|\notag\\
 &\quad\le
 \begin{cases}
 C\big((1+t)^{-\frac{1}{12}}+(1+t)^{-\frac{1}{6}}\big)
       \le C(1+t)^{-\frac{1}{12}},
               &\ell=\eps_0(1+t)^{\frac{1}{3}},\\
 C\big((1+t)^{-\frac{1}{6}}+(1+t)^{-\frac{1}{3}}\big)
       \le C(1+t)^{-\frac{1}{6}},
               &\ell=\eps_0(1+t)^{\frac{2}{3}}.
 \end{cases}
 \label{eq:away-final}
\end{align}
Combining \eqref{eq:away-final} and \eqref{eq:near-final} gives
\eqref{eq:global-rate} and \eqref{eq:away-contact}. For bounded
times, increase the constants using \eqref{eq:energy} and
\eqref{eq:shift-bound}.

\subsection{The self-similar rarefaction wave}

Let $w^r(\xi)=\max\{w_-,\min\{\xi,w_*\}\}$. For
$x=y+(1+t)w_0(y)$, we have $w(1+t,x)=w_0(y)$.
The function $w^r$ is nondecreasing and $1$-Lipschitz, so
$w^r(\frac{x}{1+t})$ lies between $w_0(y)$ and $w_*$ when $y\ge0$,
and between $w_-$ and $w_0(y)$ when $y\le0$. Hence
\begin{align*}
 \left|w(1+t,x)-w^r\left(\frac{x}{1+t}\right)\right|
 &=\left|w_0(y)-w^r\left(w_0(y)+\frac{y}{1+t}\right)\right|\\
 &\le \min\left\{\frac{|y|}{1+t},\,
              \frac{w_*-w_-}{2}(1-|\tanh y|)\right\}\\
 &\le \frac{\log(2+t)}{1+t}+\frac{C\delta_R}{(2+t)^2}.
\end{align*}
For the last bound, use the first term of the minimum for
$|y|\le\log(2+t)$ and the second for $|y|>\log(2+t)$.

Where $w^r(\frac{x}{1+t})$ and $w^r(\frac{x}{t})$ differ, $|x|\le C(1+t)$.
By Section~2.1, the state on the $1$-rarefaction curve is a
Lipschitz function of the characteristic speed. Thus, for $t\ge1$,
\begin{align*}
 &\|\bar q^{\,r}(t)-\bar q(t)\|_\infty
       =\|q^r(\frac{\cdot}{t})-q^R(t)\|_\infty\\
 &\quad\le C\left(\|w(1+t,\cdot)-w^r(\frac{\cdot}{1+t})\|_\infty
           +\|w^r(\frac{\cdot}{1+t})-w^r(\frac{\cdot}{t})\|_\infty\right)\\
 &\quad\le C\left(\frac{\log(2+t)}{1+t}+\frac{\delta_R}{(2+t)^2}
                  +\sup_{|x|\le C(1+t)}\frac{|x|}{t(1+t)}\right)\\
 &\quad\le C\frac{1+\log(2+t)}{1+t}.
\end{align*}
Since $\frac{1+\log(2+t)}{1+t}\le C(1+t)^{-\frac{1}{4}}$ for $t\ge1$,
the triangle inequality and \eqref{eq:global-rate}--\eqref{eq:away-contact}
give the asserted bounds for $\bar q^{\,r}$.
This completes the proof of Theorem~\ref{thm:decay}.


\begin{thebibliography}{99}
\bibitem{KM}
S.~Kawashima and A.~Matsumura,
Asymptotic stability of traveling wave solutions of systems for one-dimensional
gas motion,
\emph{Comm. Math. Phys.} \textbf{101} (1985), 97--127.

\bibitem{HLM}
F.~Huang, J.~Li, and A.~Matsumura,
Asymptotic stability of combination of viscous contact wave with
rarefaction waves for one-dimensional compressible Navier--Stokes
system,
\emph{Arch. Ration. Mech. Anal.} \textbf{197} (2010), 89--116.

\bibitem{KVWbaro}
M.-J.~Kang, A.~F.~Vasseur, and Y.~Wang,
Time-asymptotic stability of composite waves of viscous shock and
rarefaction for barotropic Navier--Stokes equations,
\emph{Adv. Math.} \textbf{419} (2023), Paper No.~108963.

\bibitem{KVW}
M.-J.~Kang, A.~F.~Vasseur, and Y.~Wang,
Time-asymptotic stability of generic Riemann solutions for compressible
Navier--Stokes--Fourier equations,
\emph{Arch. Ration. Mech. Anal.} \textbf{249} (2025), Paper No.~42.
\href{https://doi.org/10.1007/s00205-025-02116-w}{doi:10.1007/s00205-025-02116-w}.
\href{https://arxiv.org/abs/2306.05604}{arXiv:2306.05604}.

\bibitem{LX}
T.-P.~Liu and Z.~Xin,
Nonlinear stability of rarefaction waves for compressible Navier--Stokes
equations,
\emph{Comm. Math. Phys.} \textbf{118} (1988), 451--465.

\bibitem{HHLX}
M.~He, F.~Huang, L.~Liu, and L.~Xu,
Time-decay rates toward rarefaction waves for the one-dimensional
isentropic compressible Navier--Stokes equations,
\href{https://arxiv.org/abs/2611.01956}{arXiv:2611.01956} (2026).
\bibitem{MN}
A.~Matsumura and K.~Nishihara,
Asymptotics toward the rarefaction waves of the solutions of a
one-dimensional model system for compressible viscous gas,
\emph{Japan J. Appl. Math.} \textbf{3} (1986), 1--13.

\bibitem{HMX}
F.~Huang, A.~Matsumura, and Z.-P.~Xin,
Stability of contact discontinuities for the 1-D compressible
Navier--Stokes equations,
\emph{Arch. Ration. Mech. Anal.} \textbf{179} (2006), 55--77.

\bibitem{HM}
F.~Huang and A.~Matsumura,
Stability of a composite wave of two viscous shock waves for the
full compressible Navier--Stokes equation,
\emph{Comm. Math. Phys.} \textbf{289} (2009), 841--861.

\bibitem{MZ}
C.~Mascia and K.~Zumbrun,
Stability of small-amplitude shock profiles of symmetric hyperbolic-parabolic
systems,
\emph{Comm. Pure Appl. Math.} \textbf{57} (2004), 841--876.

\bibitem{LWX}
L.~Liu, S.~Wang, and L.~Xu,
Optimal decay rates to the contact wave for 1-D compressible
Navier--Stokes equations,
\emph{J. Math. Pures Appl.} \textbf{214} (2026), Paper No.~103948.

\bibitem{single}
L.~Chen and Y.~Wang,
Time-decay rates toward viscous shock waves for the one-dimensional
compressible Navier--Stokes--Fourier system,
\href{https://arxiv.org/abs/2611.22836v1}{arXiv:2611.22836v1} (2026).

\bibitem{LiuShock}
T.-P.~Liu,
Shock waves for compressible Navier--Stokes equations are stable,
\emph{Comm. Pure Appl. Math.} \textbf{39} (1986), 565--594.
\href{https://doi.org/10.1002/cpa.3160390502}{doi:10.1002/cpa.3160390502}.

\bibitem{HL}
X.~Huang and H.~Lee,
Time-asymptotic stability of composite wave of viscous shocks and
viscous contact wave for Navier--Stokes--Fourier equations,
\href{https://arxiv.org/abs/2504.04014v4}{arXiv:2504.04014v4} (2026).

\bibitem{KV}
M.-J.~Kang and A.~F.~Vasseur,
Contraction property for large perturbations of shocks of the
barotropic Navier--Stokes system,
\emph{J. Eur. Math. Soc.} \textbf{23} (2021), 585--638.
\href{https://doi.org/10.4171/JEMS/1018}{doi:10.4171/JEMS/1018}.

\bibitem{WW}
T.~Wang and Y.~Wang,
Nonlinear stability of planar viscous shock wave to three-dimensional
compressible Navier--Stokes equations,
\emph{J. Eur. Math. Soc.} \textbf{28} (2026), 505--557.
\href{https://doi.org/10.4171/JEMS/1486}{doi:10.4171/JEMS/1486}.

\bibitem{YangContact}
D.~Yang,
Decay rate to contact discontinuities for the 1-D compressible
Navier--Stokes system,
\emph{J. Differential Equations} \textbf{269} (2020), 6529--6558.
\href{https://doi.org/10.1016/j.jde.2020.05.004}{doi:10.1016/j.jde.2020.05.004}.

\bibitem{HW}
F.~Huang and T.~Wang,
Vanishing dissipation limit for non-isentropic Navier--Stokes equations
with shock data,
\emph{Sci. China Math.} \textbf{69} (2026), 2975--3006.
\href{https://doi.org/10.1007/s11425-025-2479-6}{doi:10.1007/s11425-025-2479-6}.

\bibitem{DHLX}
R.~Duan, F.~Huang, R.~Li, and L.~Xu,
Nonlinear stability and time-decay rates of planar rarefaction waves
for the full compressible Navier--Stokes--Fourier system in $\R^3$,
\href{https://arxiv.org/abs/2611.10749}{arXiv:2611.10749} (2026).

\bibitem{HX}
F.~Huang and L.~Xu,
Nonlinear stability of planar shock waves for the
Navier--Stokes--Fourier equations in $\R^3$,
\href{https://arxiv.org/abs/2611.11383}{arXiv:2611.11383} (2026).
\end{thebibliography}
\end{document}